\documentclass[12pt]{article}
\usepackage{}

\usepackage{epsfig}
\usepackage{latexsym}
\usepackage{caption}
\usepackage{amsfonts}
\usepackage{amssymb}
\usepackage{mathrsfs}
\usepackage{amsmath}
\usepackage{enumerate}
\usepackage{graphics}
\usepackage{MnSymbol}
\usepackage{float}
\usepackage{pict2e}
\usepackage{subfigure}
\usepackage[left]{lineno}

\usepackage{color}

\makeatletter

\newcommand{\Rmnum}[1]{\expandafter\@slowromancap\romannumeral #1@}

\makeatother

\begin{document}
	
	
	\newtheorem{theorem}{Theorem}
	\newtheorem{observation}[theorem]{Observation}
	\newtheorem{corollary}[theorem]{Corollary}
	\newtheorem{algorithm}[theorem]{Algorithm}
	\newtheorem{definition}[theorem]{Definition}
	\newtheorem{guess}[theorem]{Conjecture}
	\newtheorem{claim}[theorem]{Claim}
	\newtheorem{problem}[theorem]{Problem}
	\newtheorem{question}[theorem]{Question}
	\newtheorem{lemma}[theorem]{Lemma}
	\newtheorem{proposition}[theorem]{Proposition}
	\newtheorem{fact}[theorem]{Fact}

	\makeatletter
	\newcommand\figcaption{\def\@captype{figure}\caption}
	\newcommand\tabcaption{\def\@captype{table}\caption}
	\makeatother

	\newtheorem{acknowledgement}[theorem]{Acknowledgement}

	\newtheorem{axiom}[theorem]{Axiom}
	\newtheorem{case}[theorem]{Case}
	\newtheorem{conclusion}[theorem]{Conclusion}
	\newtheorem{condition}[theorem]{Condition}
	\newtheorem{conjecture}[theorem]{Conjecture}
	\newtheorem{criterion}[theorem]{Criterion}
	\newtheorem{example}[theorem]{Example}
	\newtheorem{exercise}[theorem]{Exercise}
	\newtheorem{notation}{Notation}
	\newtheorem{solution}[theorem]{Solution}
	\newtheorem{summary}[theorem]{Summary}

	\newenvironment{proof}{\noindent {\bf
			Proof.}}{\rule{3mm}{3mm}\par\medskip}
	\newcommand{\remark}{\medskip\par\noindent {\bf Remark.~~}}
	\newcommand{\pp}{{\it p.}}
	\newcommand{\de}{\em}
	\newcommand{\mad}{\rm mad}
	\newcommand{\qf}{Q({\cal F},s)}
	\newcommand{\qff}{Q({\cal F}',s)}
	\newcommand{\qfff}{Q({\cal F}'',s)}
	\newcommand{\f}{{\cal F}}
	\newcommand{\ff}{{\cal F}'}
	\newcommand{\fff}{{\cal F}''}
	\newcommand{\fs}{{\cal F},s}
	\newcommand{\s}{\mathcal{S}}
	\newcommand{\G}{\Gamma}
	\newcommand{\g}{(G_3, L_{f_3})}
	\newcommand{\wrt}{with respect to }
	\newcommand {\nk}{ Nim$_{\rm{k}} $  }
	\newcommand {\dom}{ {\rm Dom}  }
	\newcommand {\ran}{ {\rm Ran}  }
	\newcommand{\gs}{(G, \sigma)}
	\newcommand{\ch}{{\rm ch}}
	
	\newcommand{\q}{\uppercase\expandafter{\romannumeral1}}
	\newcommand{\qq}{\uppercase\expandafter{\romannumeral2}}
	\newcommand{\qqq}{\uppercase\expandafter{\romannumeral3}}
	\newcommand{\qqqq}{\uppercase\expandafter{\romannumeral4}}
	\newcommand{\qqqqq}{\uppercase\expandafter{\romannumeral5}}
	\newcommand{\qqqqqq}{\uppercase\expandafter{\romannumeral6}}
	
	\newcommand{\qed}{\hfill\rule{0.5em}{0.809em}}
	
	\newcommand{\var}{\vartriangle}
	
	\newcommand{\ERCagreement}{{\begin{minipage}{.60\textwidth}This paper is part of a project that has received funding from the European Research Council (ERC) under the European Union's Horizon 2020 research and innovation programme (grant agreement No 810115 -- {\sc Dynasnet}).
			\end{minipage}\hfill\begin{minipage}{.28\textwidth}\includegraphics[width=\textwidth]{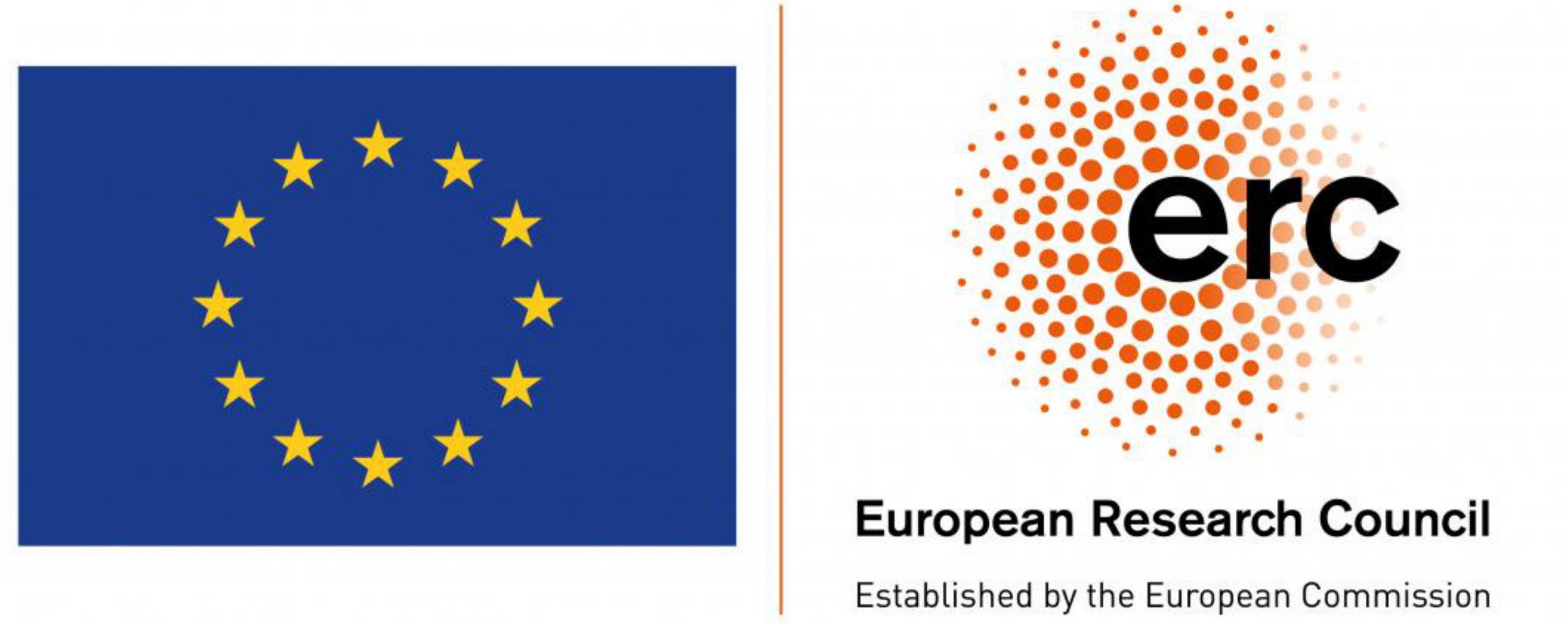}\end{minipage}\hfill}
		
		\ \ This work is partially supported by the ANR project HOSIGRA (ANR-17-CE40-0022).
		
		\ \ The  authors were also supported partially by the project 21-10775S of the Czech Science Foundation (GACR).}

	\title{{\Large \bf  Minimal asymmetric hypergraphs}\thanks\ERCagreement}
	
	\author{Yiting Jiang \thanks{Institute of Mathematics, School of Mathematical Sciences, Nanjing Normal University, Nanjing, 210023, China
			. E-mail: ytjiang@zjnu.edu.cn.} \and Jaroslav Ne\v{s}et\v{r}il\thanks{ Computer Science Institute of Charles University (IUUK and ITI), Malostransk\' e n\' am.25, 11800 Praha 1, Czech Republic. E-mail: nesetril@iuuk.mff.cuni.cz. }}
	
	\maketitle

	\begin{abstract}
		In this paper, we prove that for any $k\ge 3$, there exist infinitely many minimal asymmetric $k$-uniform hypergraphs. This is in a striking contrast to $k=2$, where it has been proved recently that there are exactly $18$ minimal asymmetric graphs.
		
		We also determine, for every $k\ge 1$, the minimum size of an asymmetric $k$-uniform hypergraph.	
		
		\bigskip
		
		\noindent {\bf Keywords:}
		asymmetric hypergraphs, $k$-uniform hypergraphs, automorphism.
	\end{abstract}
	
	\section{Introduction}
	In this paper we deal with (undirected) graphs, oriented graphs and more general hypergraphs and relational structures.
	Let us start with (undirected) graphs: An (undirected) graph $G$ is called \emph{asymmetric} if it does not have a non-identity automorphism. Any non-asymmetric graph is also called \emph{symmetric} graph. A graph $G$ is called \emph{minimal asymmetric} if $G$ is asymmetric and every non-trivial induced subgraph of $G$ is symmetric (here $G'$ is a \emph{non-trivial subgraph} of $G$ if $G'$ is a subgraph of $G$ and $1<|V(G')|<|V(G)|$). In this paper all graphs are finite.
	
	It is a folklore result that most graphs are asymmetric. In fact, as shown by Erd\H{o}s and R\'{e}nyi \cite{AGER} most graphs on large sets are asymmetric in a very strong sense. 
	The paper \cite{AGER} contains many extremal results (and problems), which motivated further research on extremal properties of asymmetric graphs, see e.g. \cite{NS}, \cite{SS}.
	This has been also studied in the context of the reconstruction conjecture \cite{CTAT}, \cite{VM}.
	
	The second author bravely conjectured a long time ago that there are only finitely many minimal asymmetric graphs, see e.g. \cite{BOOK}. Partial results were given in \cite{NS92}, \cite{Sab91}, \cite{Woj96} and recently this conjecture has been confirmed by Pascal Schweitzer and Patrick Schweitzer \cite{MSH} (the list of $18$ minimal asymmetric graphs has been isolated already in \cite{NS92}):
	
	\begin{theorem}\cite{MSH}
		\label{thm1}
		There are exactly 18 minimal asymmetric undirected graphs up to isomorphism.
	\end{theorem}
	
	\medskip
	In this paper, we consider analogous questions for \emph{$k$-graphs} (or $k$-uniform hypergraphs), i.e. pairs $(X,\mathscr{M})$ where $\mathscr{M}\subseteq {\binom{X}{k}}=\{A\subseteq X;|A|=k\}$. 
	Induced subhypergraphs, asymmetric hypergraphs and minimal asymmetric hypergraphs are defined analogously as for graphs. 
	
	We prove two results related to minimal asymmetric $k$-graphs.
	
	Denote by $n(k)$ the minimum number of vertices of an asymmetric $k$-graph.
	
	\begin{theorem}
		\label{main1}
		$n(2)=6$, $n(3)=6$, $n(k)=k+2$ for $k\ge 4$.	
	\end{theorem}

	Our second result disproves analogous minimality conjecture (i.e. a result analogous to Theorem \ref{thm1}) for $k$-graphs.
	
	\begin{theorem}
		\label{main2}
		For every integer $k\ge 3$, there exist infinitely many $k$-graphs that are minimal asymmetric.
	\end{theorem}
	
	In fact we prove the following stronger statement.
	
	\begin{theorem}
		\label{main3}
		For every integer $k\ge 3$, there exist infinitely many $k$-graphs $(X, \mathscr{M})$ such that 
		\begin{itemize}
			\item[1.] $(X, \mathscr{M})$ is asymmetric. 
			\item[2.] If $(X', \mathscr{M}')$ is a non-trivial sub-$k$-graph of $(X, \mathscr{M})$ with at least two vertices, then $(X', \mathscr{M}')$ is symmetric.
		\end{itemize}
	\end{theorem}
	
	We call $k$-graphs that satisfy the two  above properties \emph{strongly minimal asymmetric}. So strongly minimal asymmetric $k$-graphs do not contain any non-trivial (not necessarily induced) asymmetric sub-$k$-graph. Note that some of the minimal asymmetric graphs fail to be strongly minimal. For instance, as depicted in Figure \ref{fig1}, the graph $X_2$ is minimal asymmetric but not strongly minimal asymmetric, since $X_1$ is a minimal asymmetric subgraph of $X_2$.
	
	An {\emph{involution}} of a graph $G$ is any non-identity automorphism $\phi$ for which $\phi\circ\phi$ is an identity. It was proved in \cite{MSH} that all minimal asymmetric graphs are in fact minimal involution-free graphs.	
	However, it is not the case for $k$-graphs: 
	there are $k$-graphs that are (strongly) minimal asymmetric but not minimal involution-free (see examples after the proof of Theorem \ref{main3} in Section \ref{sec3}). We prove the following form of Theorem \ref{main3} relating minimal asymmetric $k$-graphs for involutions.	
	
	\begin{theorem}
		\label{main4} 
		For every integer $k\ge 6$, there exist infinitely many $k$-graphs $(X, \mathscr{M})$ such that 
		\begin{itemize}
			\item[1.] $(X, \mathscr{M})$ is asymmetric. 
			\item[2.] If $(X', \mathscr{M}')$ is a sub-$k$-graph of $(X, \mathscr{M})$ with at least two vertices, then $(X', \mathscr{M}')$ has an involution.
		\end{itemize}
	\end{theorem}
	
	\medskip
	
	Theorem \ref{main3} and Theorem \ref{main4} are proved by constructing a sequence of strongly minimal asymmetric $k$-graphs. We have two different constructions of increasing strength. In Section \ref{sec3} we give a construction with all vertex degrees bounded by $3$. A stronger construction which yields minimal asymmetric $k$-graphs ($k\ge 6$) with respect to involutions is given in the proof of Theorem \ref{main4} in Section \ref{sec4}.
	In Section \ref{sec:open} we consider minimal asymmetric relations and their multiplicities and conclude with several open problems.
	
	\section{The proof of Theorem \ref{main1}}
	
	\begin{lemma}
		\label{lem}
		For $k\ge 3$, we have $n(k)\ge k+2$. 
	\end{lemma}
	
	\begin{proof}
		Assume that there exists an asymmetric $k$-graph $(X, \mathscr{M})$ with $|X|=k+1$.
		If for each vertex $u\in X$, there is a hyperedge $M\in\mathscr{M}$ such that $u\notin M$, then $\mathscr{M}={\binom{X}{k}}$, which  is symmetric.
		Otherwise there exists $u,v\in X$ such that $\{u,v\}\subset M$  for every edge $M\in\mathscr{M}$, or there exist $u',v'\in X$ and $M_1,M_2\in\mathscr{M}$ such that $u'\notin M_1$ and $v'\notin M_2$.
		In the former case, there is an automorphism $\phi$ of $(X, \mathscr{M})$ such that $\phi(u)=v$ and $\phi(v)=u$. In the latter case there is an automorphism $\phi$ of $(X,\mathscr{M})$ such that $\phi(u')=v'$ and $\phi(v')=u'$. In either case we have a contradiction. 
	\end{proof}
	
	For a $k$-graph $G=(X,\mathscr{M})$, the \emph{set-complement} of $G$ is defined as a $(|X|-k)$-graph $\bar{G}=(X,\bar{\mathscr{M}})=(X,\{X-M|M\in\mathscr{M}\})$.
	Denote by $Aut(G)$ the set of all the automorphisms of $G$ and thus we have $Aut(G)=Aut(\bar{G})$. 
	We define the degree of a vertex $v$ in a $k$-graph $G$ as 
	$d_G(v)=|\{M\in\mathscr{M};v\in M\}|$. 
	
	\begin{lemma}
		\label{lem1}
		For $k\ge 4$, we have $n(k)=k+2$.
	\end{lemma}
	

	\begin{figure}[!htp]
		\begin{minipage}{0.32\linewidth}
			\centering
			\includegraphics[scale=0.5]{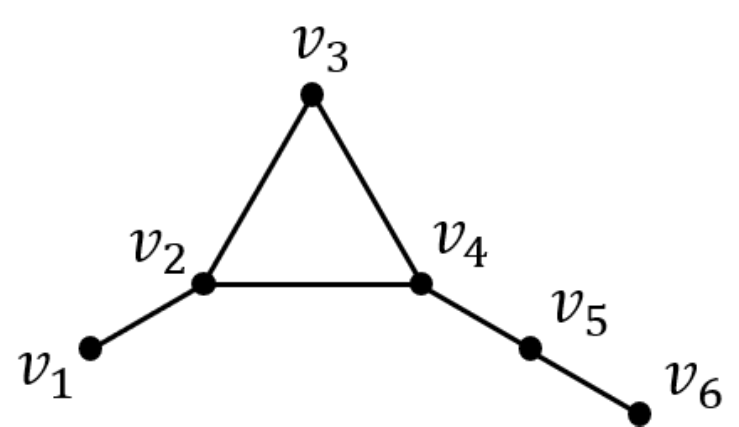}
			\caption*{$X_1$}
		\end{minipage}
		\begin{minipage}{0.32\linewidth}
			\centering
			\includegraphics[scale=0.5]{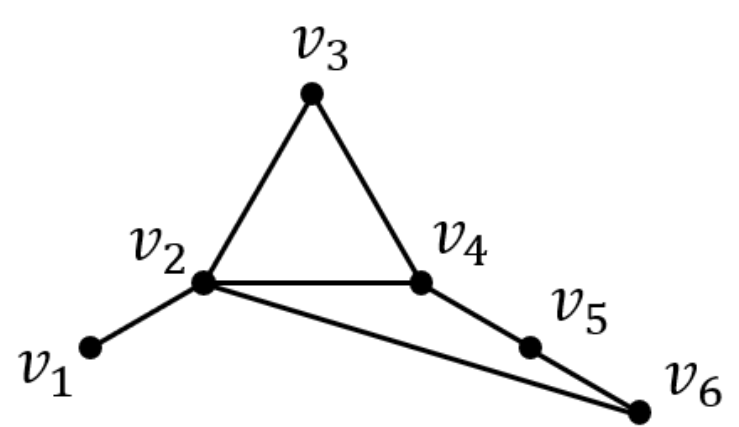}
			\caption*{$X_2$}  		
		\end{minipage}
		\begin{minipage}{0.32\linewidth}
			\centering
			\includegraphics[scale=0.5]{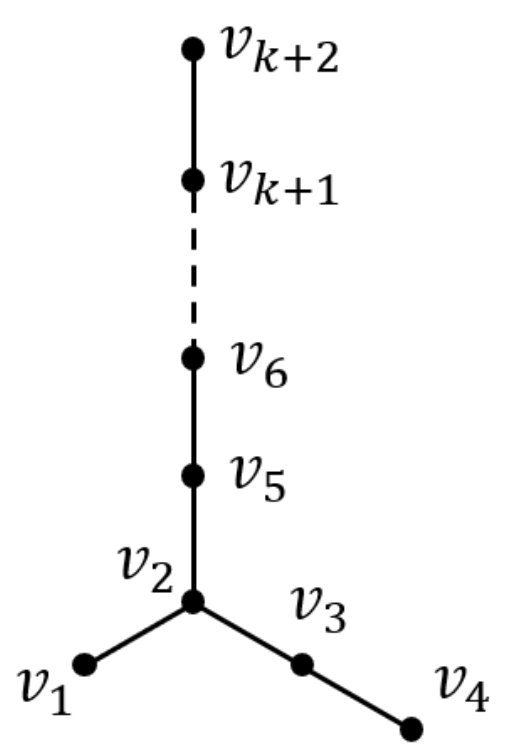}
			\caption*{$T_{k+2}$}  		
		\end{minipage}
		\caption{}
		\label{fig1}
	\end{figure}
	
	
	\begin{proof}
		First, we construct an asymmetric $2$-graph $(X,\mathscr{M})$ with $|X|=k+2$ for each $k\ge 4$.
		Examples of such graphs $X_1$ and $T_{k+2}$ are depicted in Figure \ref{fig1}.

		For $k=4$, take the set-complement of $X_1$.
		For every $k\ge 5$, take the set-complement of $T_{k+2}$.  
		It is easy to see that $X_1$ and $T_{k+2}$ ($k\ge 5$) are asymmetric.
		Thus set-complements $\bar{X_1}$ and $\bar{T_{k+2}}$ ($k\ge 5$) are also asymmetric $k$-graphs. 
		
		Each of the set-complements $\bar{X_1}$ and $\bar{T_{k+2}}$ has $k+2$ vertices. Thus a non-trivial subgraph of each of them is symmetric, by Lemma \ref{lem}.
	\end{proof}

	\begin{lemma}
		\label{lem2}
		For $k=3$, we have $n(3)=6$.
	\end{lemma}
	
	
	\begin{figure}[!htp]
		\begin{center}
			\includegraphics[scale=0.5]{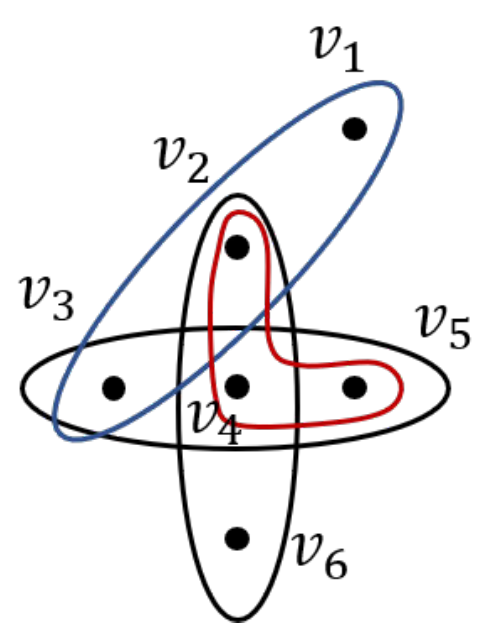}
		\end{center}
		\caption{An asymmetric $3$-graph with $|X|=6$}
		\label{fig2}
	\end{figure}


	\begin{proof}
		For $n(3)\le 6$, consider the following $3$-graph $G=(X, \mathscr{M})$ depicted on Figure \ref{fig2}, $X=\{v_1,v_2,v_3,v_4,v_5,v_6\}$, $\mathscr{M}=\{\{v_1,v_2,v_3\},\{v_2,v_4,v_5\},\{v_2,v_4,v_6\},\{v_3,v_4,v_5\}\}$.
		Observe that $d_G(v_1)=d_G(v_6)=1$, $d_G(v_2)=d_G(v_4)=3$, $d_G(v_3)=d_G(v_5)=2$. It is not difficult to see $G$ is asymmetric.
		
					Now we shall prove that $n(3)\ge 6$.
	
	Assume that there exists an asymmetric $3$-graph $H=(X, \mathscr{M})$ with $|X|=5$. Let $X=\{v_1,v_2,v_3,v_4,v_5\}$. Without loss of generality, let $M=\{v_1,v_2,v_3\}\in\mathscr{M}$. Then there exists an edge $M\in \mathscr{M}$ such that $v_4\in M$ and $v_5\notin M$, or $v_4\notin M$ and $v_5\in M$.
	
	$k$-graph $H$ is asymmetric if and only if $(X,{\binom X 3}-\mathscr{M})$ is asymmetric. Thus we can sufficiently consider that $|\mathscr{M}|\le{\frac{{\binom 5 3}}{2}}=5$. If $d_H(v_4)=d_H(v_5)$, which means both of $v_4$ and $v_5$ have degree $1$ or $2$, then there exists an automorphism $\phi$ of $H$ such that $\phi(v_4)=v_5$. Assume that $d_H(v_4)>d_H(v_5)$.
	
	
	\begin{figure}[!htp]
			\centering
		\begin{minipage}{0.3\linewidth}
			\centering
			\includegraphics[scale=0.5]{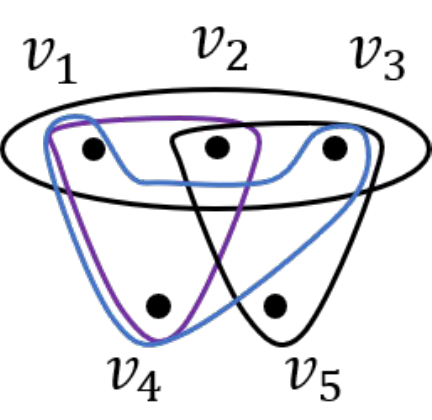}
			\caption*{(a)}
		\end{minipage}
		\begin{minipage}{0.3\linewidth}
			\centering
			\includegraphics[scale=0.5]{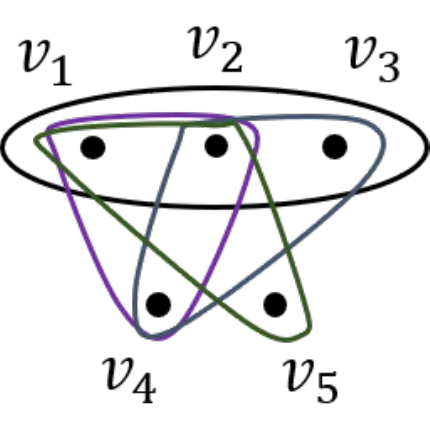}
			\caption*{(b)}  		
		\end{minipage}
		\begin{minipage}{0.3\linewidth}
			\centering
			\includegraphics[scale=0.5]{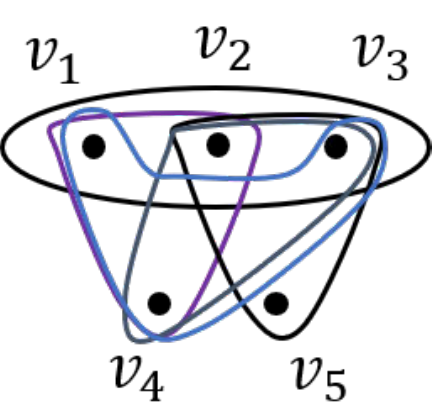}
			\caption*{(c)}
		\end{minipage}
		\caption{}
	\end{figure}

	
	{\bfseries{Case 1.}} There is no edge $M\in\mathscr{M}$ such that $\{v_4,v_5\}\subseteq M$
	
	It is sufficient to consider two subcases: $d_H(v_4)=2$ and $d_H(v_5)=1$, or $d_H(v_4)=3$ and $d_H(v_5)=1$.
	In the first subcase, up to isomorphism, we obtain two different graphs as Figure 3(a) and 3(b) shown. There exists an automorphism $\phi$ of $H$ such that $\phi(v_2)=v_3$ in (a) (resp. $\phi(v_3)=v_4$ in (b)). In the second subcase, there is only one possible graph as Figure 3(c) shown. Observe that there exists an automorphism $\phi$ of $H$ such that $\phi(v_2)=v_3$ or $\phi(v_1)=v_4$.
	
	\medskip
	
	{\bfseries{Case 2.}} There exists $M\in\mathscr{M}$ such that $\{v_4,v_5\}\subseteq M$
	
	Let $\mathbb{M}=\{M\in\mathscr{M}; \{v_4,v_5\}\subseteq M\}$, then by symmetric, $|\mathbb{M}|\ne 3$. Since $d_H(v_4)>d_H(v_5)$ and $|\mathscr{M}|\le 5$, the graphs in this case we need to consider can be divide as follow:
	\begin{itemize}
		\item[1)] $|\mathscr{M}-\mathbb{M}|=2$, as Figure 4(a), 4(b) and 4(c) shown.
		\item[2)] $|\mathscr{M}-\mathbb{M}|=3$, as Figure 4(d), 4(e) and 4(f) shown.
		\item[3)] $|\mathscr{M}-\mathbb{M}|=4$, as Figure 4(g), 4(h), 4(i), 4(j), 4(k) and 4(l) shown.
	\end{itemize}
	It is easily to observe that there is an automorphism $\phi$ such that $\phi(v_1)=v_4$ and $\phi(v_3)=v_5$ in Figure 4(a), $\phi(v_1)=v_5$ and $\phi(v_2)=v_4$ in Figure 4(b), $\phi(v_1)=v_2$ in Figure 4(c), $\phi(v_1)=v_3$ in Figure 4(d), $\phi(v_1)=v_4$ in Figure 4(e), $\phi(v_2)=v_4$ and $\phi(v_3)=v_5$ in Figure 4(f), $\phi(v_1)=v_3$ or $\phi(v_2)=v_4$ in Figure 4(g), $\phi(v_1)=v_3$ in Figure 4(h), $\phi(v_3)=v_4$ in Figure 4(i), $\phi(v_1)=v_4$ or $\phi(v_3)=v_5$ in Figure 4(j), $\phi(v_1)=v_2$ and $\phi(v_3)=v_5$ in Figure 4(k), $\phi(v_3)=v_4$ in Figure 4(l). In each case, we obtain a contradiction.
\end{proof}


\begin{figure}[!htp]
			\centering
	\begin{minipage}{0.3\linewidth}
		\centering
		\includegraphics[scale=0.5]{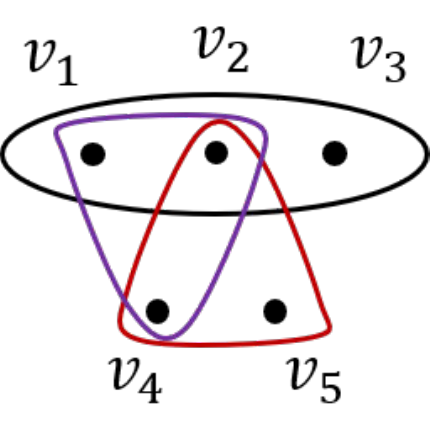}
		\caption*{(a)}
	\end{minipage}
	\begin{minipage}{0.3\linewidth}
		\centering
		\includegraphics[scale=0.5]{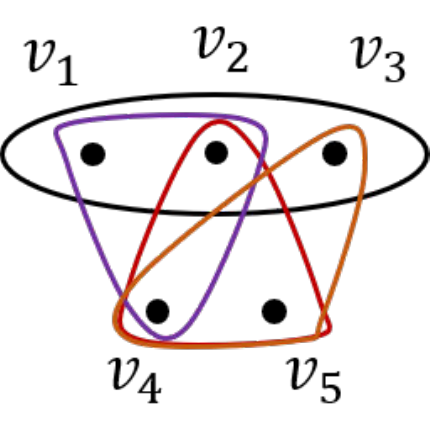}
		\caption*{(b)}  		
	\end{minipage}
	\begin{minipage}{0.3\linewidth}
		\centering
		\includegraphics[scale=0.5]{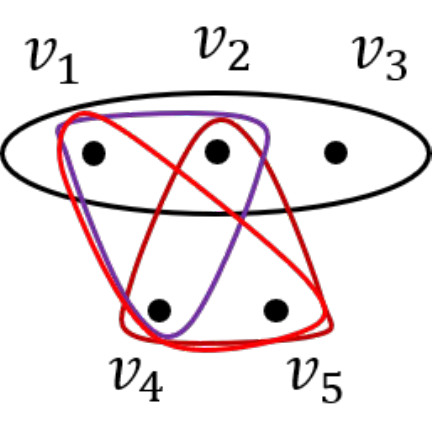}
		\caption*{(c)}
	\end{minipage}	
	\begin{minipage}{0.3\linewidth}
		\centering
		\includegraphics[scale=0.5]{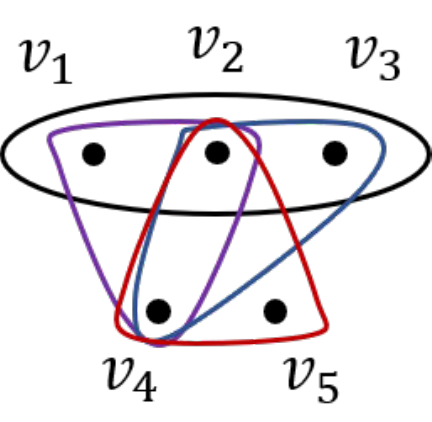}
		\caption*{(d)}
	\end{minipage}
	\begin{minipage}{0.3\linewidth}
		\centering
		\includegraphics[scale=0.5]{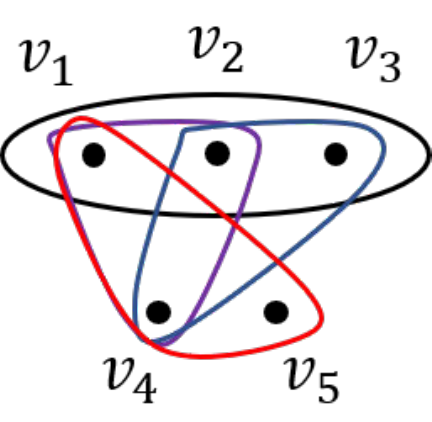}
		\caption*{(e)}  		
	\end{minipage}
	\begin{minipage}{0.3\linewidth}
		\centering
		\includegraphics[scale=0.5]{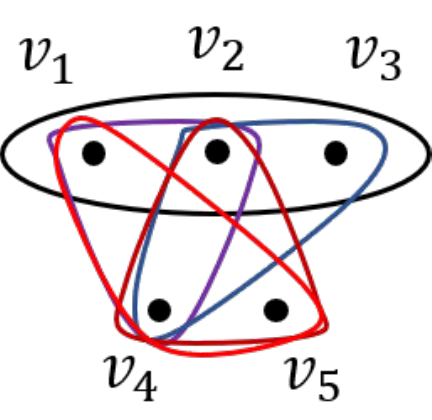}
		\caption*{(f)}
	\end{minipage}	
	\begin{minipage}{0.3\linewidth}
		\centering
		\includegraphics[scale=0.5]{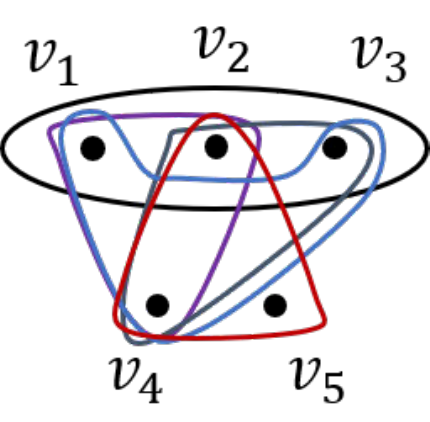}
		\caption*{(g)}
	\end{minipage}
	\begin{minipage}{0.3\linewidth}
		\centering
		\includegraphics[scale=0.5]{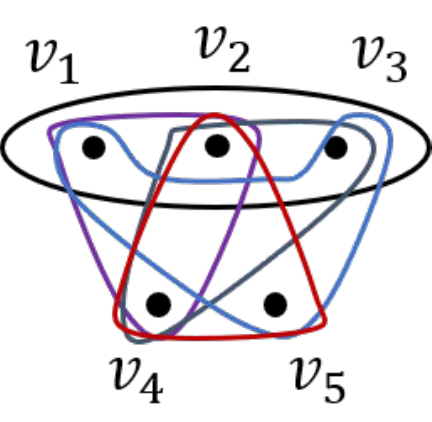}
		\caption*{(h)}  		
	\end{minipage}
	\begin{minipage}{0.3\linewidth}
		\centering
		\includegraphics[scale=0.5]{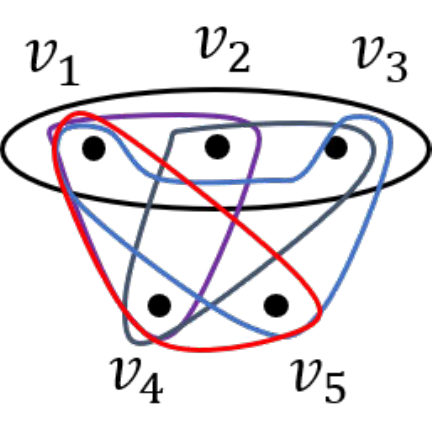}
		\caption*{(i)}  		
	\end{minipage}
	\begin{minipage}{0.3\linewidth}
		\centering
		\includegraphics[scale=0.5]{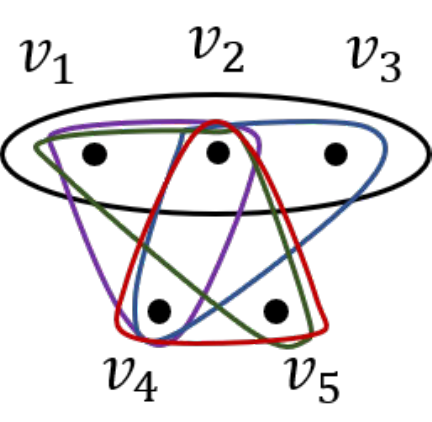}
		\caption*{(j)}  		
	\end{minipage}
	\begin{minipage}{0.3\linewidth}
		\centering
		\includegraphics[scale=0.5]{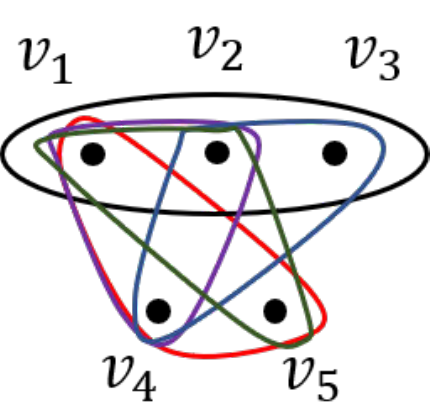}
		\caption*{(k)}  		
	\end{minipage}
	\begin{minipage}{0.3\linewidth}
		\centering
		\includegraphics[scale=0.5]{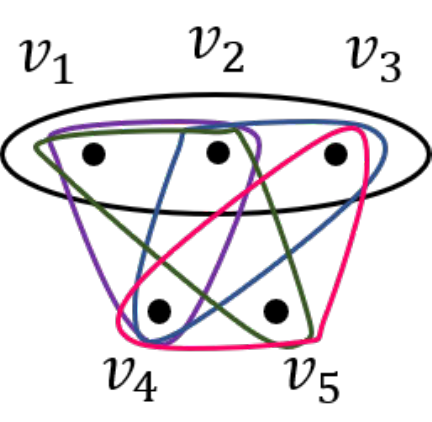}
		\caption*{(l)}  		
	\end{minipage}	
	\caption{}
\end{figure}

	\section{Minimal asymmetric $k$-graphs}
	
	In this section, we give proofs of Theorem \ref{main3} and Theorem \ref{main4}.

	\subsection{Proof of Theorem \ref{main3}}
	\label{sec3}
	
	We define the following $k$-graphs for $k\ge 3$, $t\ge k-2$. (Note that for each positive integer $p$, we denote by $[p]$ the set $\{0,1,2,\dots,p-1\}$.)
	
	\medskip
	$G_{k,t}=(X_{k,t},\mathscr{E}_{k,t})$, 
	
	$X_{k,t}=\{u_i;i\in [tk]\}\cup\{v^{j}_i;i\in [tk], j\in [k-2]\}\}$,
	
	$\mathscr{E}_{k,t}=\{E_i;i\in [tk]\}\cup\{E_{i,j};j\in\{1,2,\dots,k-3\},i=j+sk-1,s\in[t]\}$, \\where $E_i=\{v^0_i,u_i,v^1_i,v^2_i,\dots,v^{k-3}_i,v^0_{i+1}\}$, $E_{i,j}=\{v^j_i,v^j_{i+1},\dots,v^j_{i+k-1}\}$, 
	and using addition modulo $tk$.
	
	\medskip
	$G^{\circ}_{k,t}=\{X_{k,t}\cup\{x\}, \mathscr{E}_{k,t}\cup\{E^0\}\}$, where $E^0=\{v^0_0,u_0,v^1_0,v^2_0,\dots,v^{k-3}_0,x\}$.

\begin{figure}[!htp]
		\begin{minipage}{0.48\linewidth}
			\centering
			\includegraphics[scale=0.3]{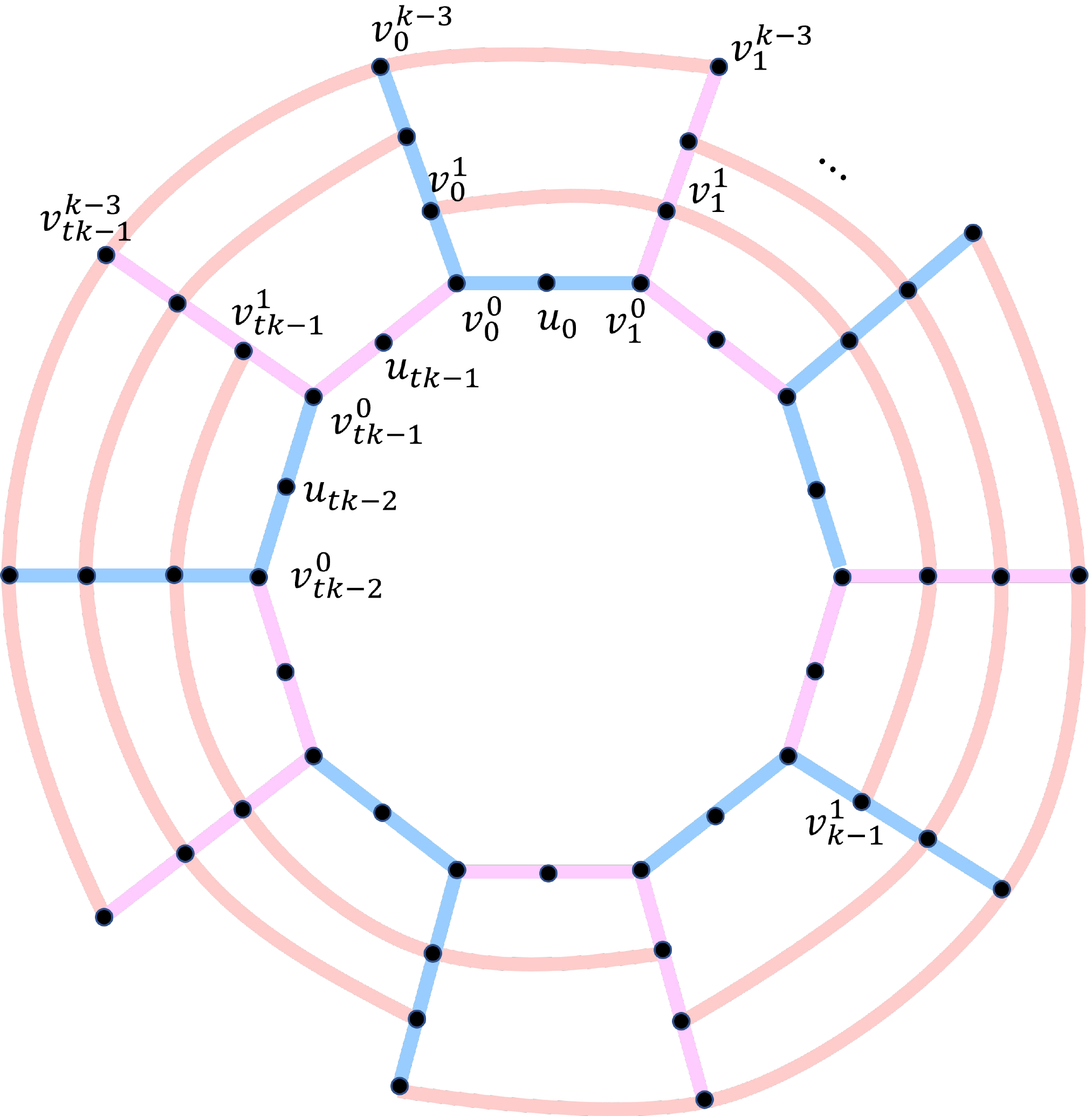}
			\caption*{$G_{k,t}$}
		\end{minipage}
		\begin{minipage}{0.48\linewidth}
			\centering
			\includegraphics[scale=0.3]{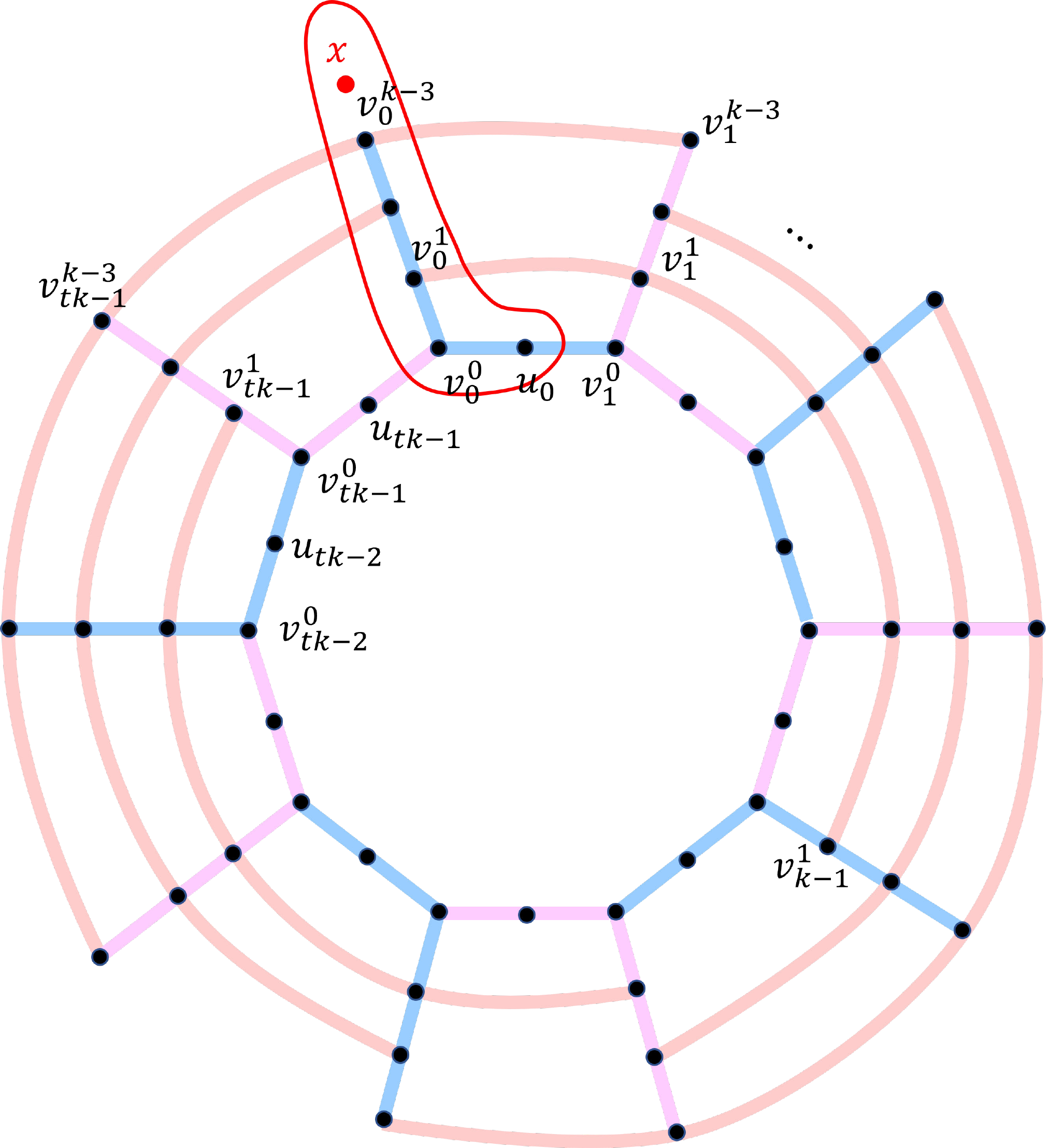}
			\caption*{$G^{\circ}_{k,t}$}  		
		\end{minipage}
		\caption{}
		\label{fig3}
	\end{figure}

	\medskip
	The graphs $G_{k,t}$ and $G^{\circ}_{k,t}$ is schematically depicted on Figure \ref{fig3}.
	
	
	\medskip
	The proof of Theorem \ref{main3} follows from the following two lemmas.
	
	\begin{lemma}
		\label{lem7}
		\begin{itemize}
			\item[1)] The graph $G_{k,t}$ is symmetric and every non-identity automorphism $\phi$ of $G_{k,t}$ satisfies one of the following properties.
   \begin{itemize}
       \item There exists a positive integer $c\ne 0 \mod tk$ such that for every $i\in [tk]$, $j=(i+c)$, $\phi(E_i)=E_j$ (i.e. for each vertex $v\in E_i$, $\phi(v)\in E_j$);
       \item There exists an $i\in[tk]$ such that $\phi(E_i)=E_{i+1}$;
       \item There exists an $i\in[tk]$ such that $\phi(E_i)=E_{i+2}$.
   \end{itemize}
			\item[2)] The only automorphism of $G_{k,t}$ which leaves the set $E_0\setminus\{v_1\}$ invariant (i.e. for each vertex $v\in E_0\setminus\{v_1\}$, $\phi(v)\in E_0\setminus\{v_1\}$) is the identity. 
			\item[3)] Every non-trivial subgraph of $G_{k,t}$ containing the vertices in $E_0$ has a non-identity automorphism $\phi$ which leaves the set $E_0$ invariant.
		\end{itemize}
	\end{lemma}
	
	\begin{proof}
		The first property can be seen to hold by considering the degrees of the vertices in $G_{k,t}$. The second property follows easily from this.
		
		To prove the third property, let $G$ be a non-trivial subgraph of $G_{k,t}$ containing the vertices in $E_0$ and let $s$ be the maximal index such that $G$ contains the edges $E_0$, $E_1$, $\dots$, $E_s$. Suppose first that $s\ne tk-1$. Since $E_{s+1}$ is not in $G$, the vertices $u_s$ and $v_{s+1}$ are of degree one. 
		The automorphism $\phi$ of $G$ which interchanges $u_s$ and $v_{s+1}$ and leaves all the other vertices fixed is a (non-identity) involution. If $s=kt-1$, there is an edge $E_l$, $l\in\{i,i+1,\dots,i+k-1\}$, the vertices $u_l$ and $v^j_l$ have degree one. So here also there is a (non-identity) involution of $G$ that interchanges $u_l$ and $v^j_l$ and leaves all other vertices fixed, in particular leaving $E_0$ invariant.
	\end{proof}
	
	\begin{lemma}
		\begin{itemize}
			\item[1)] The graph $G^{\circ}_{k,t}$ is asymmetric. 
			\item[2)] Every non-trivial subgraph of $G^{\circ}_{k,t}$ has a non-identity automorphism.
		\end{itemize}
	\end{lemma}	
	
	\begin{proof}
		To prove the first property, we first suppose that $\phi$ is an automorphism of $G^{\circ}_{k,t}$. We can see that the edges $E^0$ and $E_0$ are invariant under $\phi$ by considering the degrees of the vertices in $G^{\circ}_{k,t}$. Since $G_{k,t}$ is a subgraph of $G^{\circ}_{k,t}$, the automorphism $\phi'$ induced by $\phi$ on $G_{k,t}$ leaves the $E_0\setminus\{v_1\}$ invariant. By Lemma \ref{lem7}, $\phi'$ is identity, thus $\phi$ is identity. Therefore, $G^{\circ}_{k,t}$ is asymmetric.

		To prove the second property, let $G$ be a non-trivial subgraph of $G^{\circ}_{k,t}$. If $G$ contains the edge $E_0$, then either $G=G_{k,t}$ or $G$ contains a non-trivial subgraph of $G_{k,t}$ containing the vertices in $E_0$. In both of the cases, according to Lemma \ref{lem7}, there is a non-identity automorphism of $G$. Suppose that $G$ does not contain the edge $E_0$. If $E^0$ is in $G$, then there is a non-identity involution of $G$ that interchanges $x$ and $u_0$ and leaving all other vertices fixed. If $E^0$ is not in $G$, then either $G$ does not contain any edge $E_i$ for all $i\in[tk]$ or there exist some $i\in [tk]\setminus\{0\}$ such that $E_i$ is an edge of $G$. In the former case, $G$ is consists of some pairwise disjoint edges, which is trivially symmetric. In the latter case, let $s$ be the minimal index such that $E_s$ is an edge of $G$. Since $E_{s-1}$ is not in $G$, there is a non-identity involution of $G$ that interchanges $v^0_s$ and $u_s$ and leaving all other vertices fixed.
	\end{proof}

	\medskip

	It is easy to observe that the $k$-graphs $G^{\circ}_{k,t}$ have vertex degrees at most three.
	However note that in this construction, some of the strongly minimal asymmetric $k$-graphs $G^{\circ}_{k,t}$ are not minimal involution-free. In fact, when $k\ge 3$, $t\ge k-2$ is odd, the sub-$k$-graph $G^{\circ}_{k,t}-x$ of $G^{\circ}_{k,t}$ is involution-free. 
	The most interesting form of Theorem \ref{main3} relates to minimal asymmetric graphs for involutions. It will be proved next.

	\subsection{Proof of Theorem \ref{main4}}
	\label{sec4}
	
	Let us recall Theorem \ref{main4}.
	
	\smallskip
	\noindent
	{\bfseries{Theorem \ref{main4}}}
	For every $k\ge 6$, there exist infinitely many $k$-graphs $(X, \mathscr{M})$ such that 
	\begin{itemize}
		\item[1.] $(X, \mathscr{M})$ is asymmetric. 
		\item[2.] If $(X', \mathscr{M}')$ is a sub-$k$-graph of $(X, \mathscr{M})$ with at least two vertices, then $(X', \mathscr{M}')$ has an involution.
	\end{itemize}
	(So we claim infinitely many strongly minimal involution-free $k$-graphs for every $k\ge 6$.)
	
	\medskip
	In the proof, we first construct the following $k$-graphs for $k\ge 4$:
	
	\medskip
	$G_k=(X_k,\mathscr{M}_k)$, $X_k=\{v_1,v_2,\dots,v_{2k-1}\}$, $\mathscr{M}_k=\{M_i=\{v_i,v_{i+1},\dots,v_{i+k-1}\};i\in\{1,2,\dots,k\}\}$. 
	
	\medskip
	$G^*_k=(X^*_k,\mathscr{M}^*_k)$, $X^*_k=X_k\cup\{x\}$, $\mathscr{M}^*_k=\mathscr{M}_k\cup\{M^*\}$, where $M^*=\{x,v_1,\dots,v_{k-2},v_{k+2}\}$.
	
	\medskip
	
\begin{figure}[!htp]
		\begin{center}
			\includegraphics[scale=0.3]{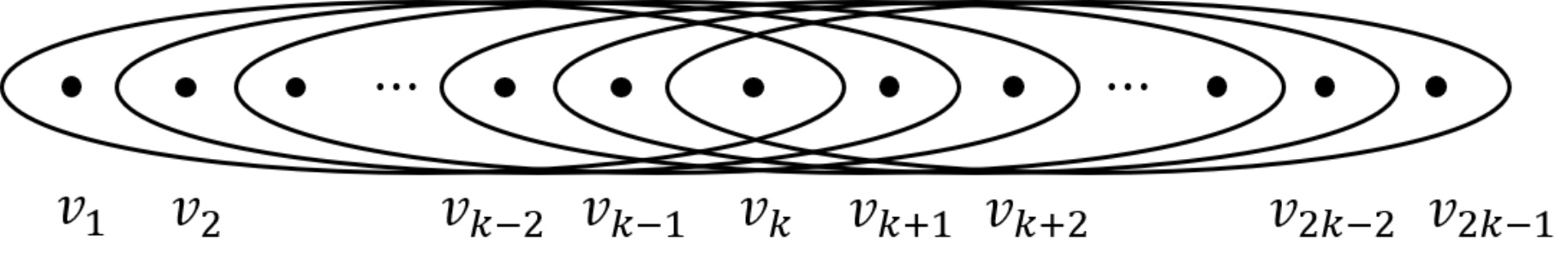}
		\end{center}
		\caption{The graph $G_k$}
		\label{fig4}
	\end{figure}

\begin{figure}[!htp]
		\begin{center}
			\includegraphics[scale=0.3]{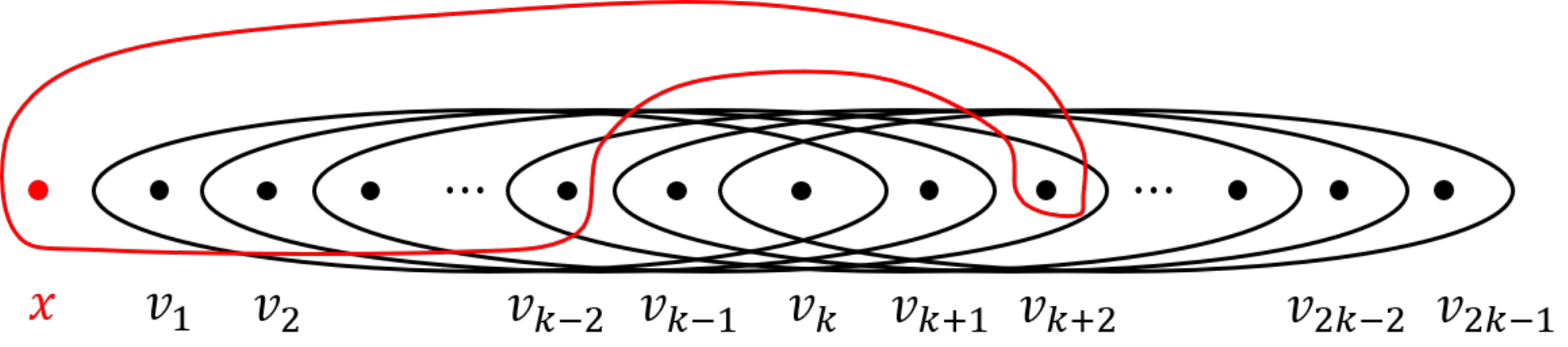}
		\end{center}
		\caption{The graph $G^*_k$}
		\label{fig5}
	\end{figure}

	These $k$-graphs are depicted on Figure \ref{fig4} and \ref{fig5}.

	
	They will be used as building blocks of our construction.
	
	\begin{lemma}
		\label{lem3}
		\begin{itemize}
			\item[1)] The $k$-graph $G_k$ is symmetric and the only non-identity automorphism $\phi$ of $G_k$ satisfies that $\phi(v_i)=v_{2k-i}$ for every $i\in\{1,2,\dots,2k-1\}$.
			\item[2)] The only automorphism of $G_k$ which leaves the set $\{v_{2k-2}, v_{2k-1}\}$ invariant (i.e. $\{\phi(v_{2k-2}), \phi(v_{2k-1})\}=\{v_{2k-2}, v_{2k-1}\}$) is the identity. 
			\item[3)] Every non-trivial sub-$k$-graph of $G_k$ containing vertices $v_{2k-2}$, $v_{2k-1}$ has an involution $\phi$ which leaves the set $\{v_{2k-2}, v_{2k-1}\}$ invariant.
			\item[4)] Every non-trivial sub-$k$-graph $G$ of $G_k$ with at least two vertices has a non-identity automorphism $\phi$, which is an involution (i.e. $\phi\circ\phi=1_{V(G)}$).
		\end{itemize}
	\end{lemma}
	
	\begin{proof}
		The first property holds by considering the degree of each vertex in $G_k$. Then also the second property follows.
		
		To prove the third one, we assume that $G$ is a non-trivial sub-$k$-graph of $G_k$ such that $G$ contains vertices $v_{2k-2}$, $v_{2k-1}$ and $j$ is the maximal index such that $G$ contains the edge $M_j=\{v_j,v_{j+1},\dots,v_{j+k-1}\}$. Let $i$ be the minimal index such that $G$ contains the edges $M_i, M_{i+1}, \dots, M_j$. Since $G$ is a nontrival sub-$k$-graph of $G_k$, we have $j<k$ and $M_{j+1}$ is not in $G$ or $i>1$ and $M_{i-1}$ is not in $G$. It implies that $v_{i+k-2}$, $v_{i+k-1}$ share the same edges $M_i$, $M_{i+1}$, $\dots$, $M_j$. If $i\notin\{k-1,k\}$ then there is an involution $\phi$ of $G$ which leaves the set $\{v_{2k-2}, v_{2k-1}\}$ invariant, $\phi(v_{i+k-2})=v_{i+k-1}$ and $\phi(v_{i+k-1})=v_{i+k-2}$. If $i\in\{k-1,k\}$ and $G$ contains an edge $M_l$ ($1\le l<i-1$), then there is an involution $\phi$ of $G$ which leaves the set $\{v_{2k-2}, v_{2k-1}\}$ invariant, $\phi(v_{i-2})=v_{i-1}$ and $\phi(v_{i-1})=v_{i-2}$, as $M_l$ contains the vertex $v_{i-1}$ and $M_{i-1}$ is not in $G$. Now the remaining case is the edge set of $G$ is contained in $\{M_{k-1},M_k\}$, which is easy to observe that there is an involution $\phi$ of $G$ which leaves the set $\{v_{2k-2}, v_{2k-1}\}$ invariant.
		
		As the proof of the last property is similar to the previous one we omit it. 
	\end{proof}
	
	\begin{lemma}
		\label{lem4}
		\begin{itemize}
			\item[1)] The $k$-graph $G^*_k$ is asymmetric. 
			\item[2)] Every non-trivial sub-$k$-graph of $G^*_k$ has an involution.
		\end{itemize}
	\end{lemma} 	
	
	\begin{proof}
		First, we prove that $G^*_k$ is asymmetric. Assume that $\phi$ is a non-identity automorphism of $G^*_k$. By considering the degrees of the vertices in the edges $M^*$ and $M_k$ we conclude that $\phi(x)=x$ and $\phi(v_{2k-1})=v_{2k-1}$ since $x$ and $v_{2k-1}$ are the only two vertices in $G^*_k$ with degree one. As $G_k$ is a sub-$k$-graph of $G^*_k$, by Lemma \ref{lem3}, we know that $\phi(v_i)= v_i$ for every $i\in\{1,2,\dots,k\}$. Thus $G^*_k$ is asymmetric. (Here one needs $k\ge 4$, which leads below to $k\ge 6$).
		
		To prove the second property of $G^*_k$,
		we assume $G$ is a non-trivial sub-$k$-graph of $G^*_k$. Then either $G$ is a sub-$k$-graph of $G_k$ or $G$ is obtained by adding the vertex $x$ and the edge $M^*=\{x,v_1,\dots,v_{k-2},v_{k+2}\}$ to a non-trivial sub-$k$-graph of $G_k$. In the former case, $G$ has an involution by Lemma \ref{lem3}. 
		In the latter case, since $G$ contains $M^*$, if there exists some $i\in\{1,2,\dots,k-2\}$ such that $M_i$ is not an edge of $G$, then $G$ has an involution $\phi$ with $\phi(x)=v_i$ and $\phi(v_i)=x$. Thus $G$ contains all of the edges $M_1, M_2,\dots, M_{k-2}$.  Let $j$ be the maximal index that $G$ contains the edges $M_1, M_2, \dots, M_j$. 
		Since $G$ is a nontrival sub-$k$-graph of $G^*_k$, we have $j<k$ and $M_{j+1}$ is not in $G$, hence $j\in\{k-2,k-1\}$. If $M_k$ is not an edge of $G$, then either $j=k-2$ or $j=k-1$ there is an involution $\phi$ such that $\phi(v_{k-1})=v_k$ and $\phi(v_k)=v_{k-1}$. Thus $G$ contains all the edges of $G^*_k$ but $M_{k-1}$. 
		So there is a (non-identity) involution of $G$ that interchanges $v_{2k-2}$ and $v_{2k-1}$ and leaves all other vertices fixed.
	\end{proof}

	For a hypergraph $G=(X,\mathscr{M})$, let $\widetilde{G}=(\widetilde{X}, \widetilde{\mathscr{M}})$ be a hypergraph with $\widetilde{X}=X\cup \bigcup\limits_{i=1}^{|\mathscr{M}|}\{a_i,b_i\}$ 
	(where $\{a_i,b_i\}\cap\{a_j,b_j\}=\emptyset$ and $\{a_i,b_i\}\cap X=\emptyset$ for any $i,j\in[|\mathscr{M}|]$)
	and $\widetilde{\mathscr{M}}=\{M_i\cup\{a_i,b_i\}; M_i\in \mathscr{M}\}$.
	
	\begin{observation}
		\label{obs1}
		For every hypergraph $G=(X,\mathscr{M})$, every automorphism of $\widetilde{G}$ which maps $X$ to $X$ is also an automorphism of $G$ and every automorphism of $G$ extends to an automorphism of $\widetilde{G}$.
	\end{observation}
	
	\begin{lemma}
		\label{lem5}
		Suppose $\phi$ is an automorphism of $\widetilde{G}_k=(\widetilde{X}_k,\widetilde{\mathscr{M}}_k)$ which leaves the set $\{v_1,v_{2k-2},v_{2k-1}\}$ invariant. Then $\phi$ restricted to $X_k$ is identity.
	\end{lemma}
	
	\begin{proof}
		Observe that the degree of each vertex in $X_k\setminus\{v_1,v_{2k-2},v_{2k-1}\}$ in $\widetilde{G}_k$ is at least $2$ while every vertex in $\widetilde{X}_k\setminus X_k$ has degree one.
		As $\phi$ is an automorphism of $\widetilde{G}_k$ which leaves the set $\{v_1,v_{2k-2},v_{2k-1}\}$ invariant, $\phi$ maps $X_k$ to $X_k$.
		By Lemma \ref{lem3} and Observation \ref{obs1}, $\phi$ restricted to $X_k$ is identity. 
	\end{proof}
	
	\begin{lemma}
		\label{lem6}
		Suppose $\phi$ is an automorphism of $\widetilde{G}^*_k=(\widetilde{X}^*_k,\widetilde{\mathscr{M}}^*_k)$ which leaves the vertices $x$ and $v_{2k-1}$ invariant. Then $\phi$ restricted to $X^*_k$ is identity.
	\end{lemma}
	
	\begin{proof}
		The proof of this lemma is very similar to the above proof of Lemma \ref{lem5}.
		Observe that the degree of each vertex in $X^*_k\setminus\{x,v_{2k-1}\}$ is at least $2$ while every vertex in $\widetilde{X}^*_k\setminus X_k$ has degree one.
		As $\phi$ is an automorphism of $\widetilde{G}^*_k$ which leaves the set $\{x,v_{2k-1}\}$ invariant, $\phi$ maps $X^*_k$ to $X^*_k$.
		By Lemma \ref{lem4} and Observation \ref{obs1}, $\phi$ restricted to $X^*_k$ is identity. 
	\end{proof}

\begin{figure}[!htp]
		\begin{center}
			\includegraphics[scale=0.25]{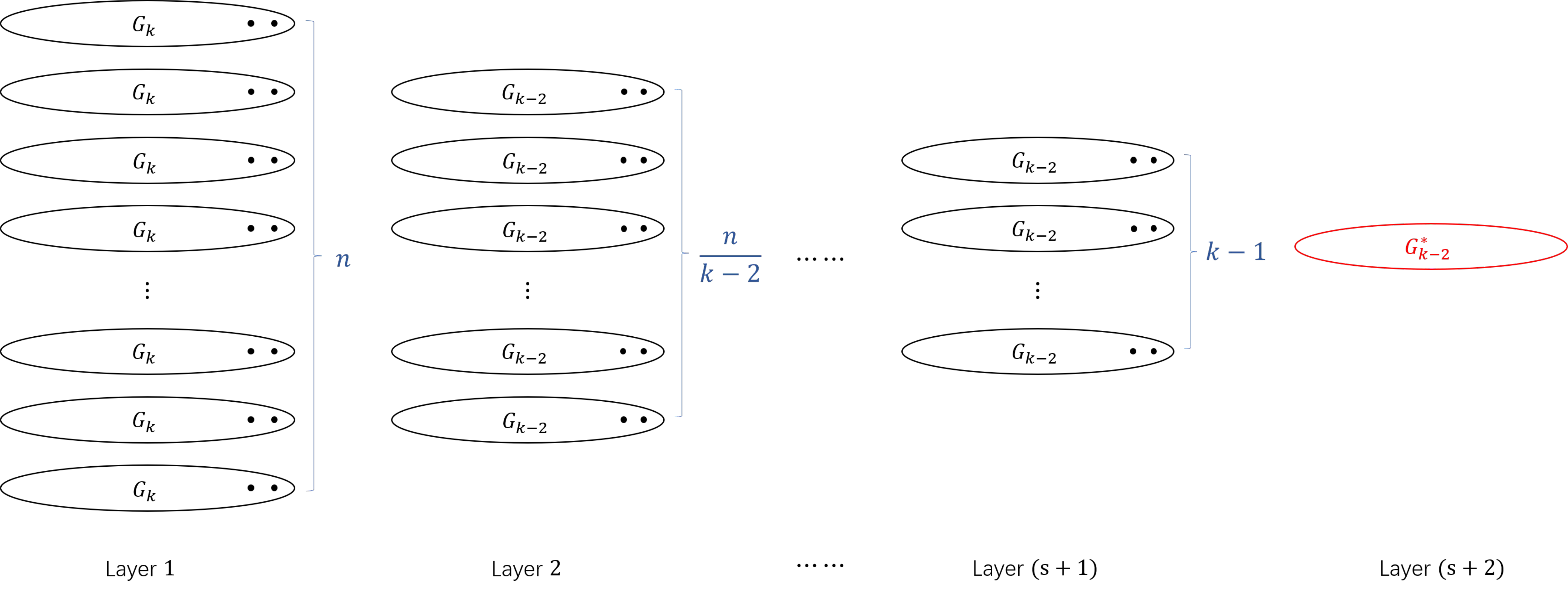}
		\end{center}
		\caption{The hypergraph $H$}
		\label{fig6}
	\end{figure}

	After all these preparations we shall,
	for each $k\ge 6$ and any non-negative integer $s$,  construct a $k$-graph $G_{k,s}=(X,\mathscr{M})$ with desired properties.
	Let $n=(k-1)(k-2)^s$.
	First, we construct a hypergraph $H=(X,\hat{\mathscr{M}})$, depicted as Figure \ref{fig6}, which is consist of $s+2$ layers as follow: 
	\begin{itemize}
		\item On layer 1, disjoint union of $n$ copies of $G_k$. 
		\item On layer 2, disjoint union of $\frac{n}{k-2}$ copies of $G_{k-2}$.
		\item On layer 3, disjoint union of $\frac{n}{(k-2)^2}$ copies of $G_{k-2}$.
		\item ...
		\item On layer $(s+1)$, disjoint union of $\frac{n}{(k-2)^s}=k-1$ copies of $G_{k-2}$.
		\item On layer $(s+2)$, one copy of $G^*_{k-2}$.
	\end{itemize}
	
	Intuitively, $G_{k,s}$ is obtained from $H$ by associating 
 to each $(k-2)$-edge in each copy of $G_{k-2}$ on layer $(i+1)$ (or $G^*_{k-2}$ on the last layer $(s+2)$) a copy of $G_{k-2}$ on layer $i$, $i\in\{1,2,\dots,s+1\}$ (or $G_k$ on layer $1$) and changing each $(k-2)$-edge into a $k$-edge by adding the last two vertices of the corresponding copy of $G_{k-2}$ (or $G_k$) to it.

	
	Formally, the $k$-graph $G_{k,s}=(X,\mathscr{M})$ can be constructed in two steps as follows. As above, set $n=(k-1)(k-2)^s$.  Consider first $n$ copies of $G_k$, ${\frac{n}{(k-2)}}+{\frac{n}{(k-2)^2}}+\cdots+(k-1)$ copies of $G_{k-2}$ and one copy of $G^*_{k-2}$ arranged into $s+2$ layers (see schematic Figure \ref{fig6}). We then have hypergraph $G^*_{k-2}$ on layer $(s+2)$. Graphs on layer $(s+1)$ are $k-1$ copies of $G_{k-2}$, which will be listed as $G(1)$, $G(2)$, $\dots$, $G(k-1)$.
	Graphs on layer $l$, $s+1\ge l\ge 1$, will be $\frac{n}{(k-2)^{l-1}}$ copies of $G_{k-2}$ (or $G_k$ when $l=1$) and they will be listed as $G(i_l,i_{l+1},\dots,i_{s+1})$, $1\le i_j\le k-2$, $j=l,l+1,\dots,s$, $1\le i_{s+1}\le k-1$.
	Then the vertices of $G_{k,s}$ are obtained from the vertices of the disjoint union of all hypergraphs $G(i_l,i_{l+1},\dots,i_{s+1})$, $1\le l\le s+1$ and $G^*_{k-2}$. All this can be made more precise at the cost of more notation. We leave this to the interested reader.
	
	Next, modify the $(k-2)$-edges to $k$-edges, which enlarge $G(i_l,i_{l+1},\dots,i_{s+1})$, $2\le l\le s+1$ and $G^*_{k-2}$ to $\widetilde{G}(i_l,i_{l+1},\dots,i_{s+1})$, $2\le l\le s+1$ and $\widetilde{G}^*_{k-2}$, as follows. We start with the last layer $s+2$.
	
	Recall that $\mathscr{M}_{k-2}=\{M_1,M_2,\dots,M_{k-2}\}$ and $\mathscr{M}^*_{k-2}=\{M^*_1,M^*_2,\dots,M^*_{k-2},M^*_{k-1}\}$.
	The edge $M^*_{i_l}$ of $G^*_{k-2}$ is enlarged by two last vertices of each hypergraph $G(i_{s+1})$ (on layer $s+1$), $1\le i_{s+1}\le k-1$. 
 
 The previous layer $l$, $2\le l\le s+1$, are treated similarly: the edge $M$ corresponding to $M_{i_{i-1}}$ in $G(i_l,\dots,i_{s+1})$ on layer $l$ is enlarged by the last two vertices of $G(i_{l-1},i_l,\dots,i_{s+1})$ on layer $l-1$, $1\le i_l\le k-2$. 
 
 This finishes the construction of the $k$-graph $G_{k,s}$. And it is easy to observe that all vertices of $k$-graphs $G_{k,s}$ have degree bounded by $k$.
	
	In the remaining of the proof, we  use $G(i_1,i_2,\dots,i_{s+1})$, $\widetilde{G}(i_l,i_{l+1},\dots,i_{s+1})$, $2\le l\le s+1$ and $\widetilde{G}^*_{k-2}$ as the corresponding sub-$k$-graphs of $G_{k,s}$. The corresponding vertex sets of $G_{k,s}$ are denoted by $V(G(i_1,i_2,\dots,i_{s+1}))$ $1\le l\le s+1$ and $V(G^*_{k-2})$.
	
	\medskip
	
	Since $s$ can be any non-negative integer, it is sufficient to prove that for each of $k$-graphs $G_{k,s}$, the properties in Theorem \ref{main4} hold.
	
	First, we prove asymmetry. To the contrary, we assume that $G_{k,s}$ has a non-identity automorphism $\phi$.

	\begin{claim}
		\label{clm1}
		$\phi(v)=v$ for every vertex $v$ on layer $(s+2)$,
	\end{claim}
	
	\begin{proof}
		Observing the degree sequence of each edge in $G_{k,s}$, the degree one belongs to four different types of degree sequences: 
  \begin{itemize}
      \item[--] the first edge in each copy $G(i_1,i_2,\dots,i_{s+1})$ of $G_k$: $(1,2,\dots,k)$;
      \item[--]  the corresponding edge $M$ of edge $M_{i_l}$ in each copy $G(i_l,i_{l+1},\dots,i_{s+1})$ of $G_{k-2}$ on layer $l$ ($2\le l\le s+1$): $(1,2,2,3,3\dots,k-2)$;
      \item[--] the two corresponding edges $M'$ and $M''$ of two different edges $M^*_1$ and $M^*_{k-1}$ of $G^*_{k-2}$ on layer $(s+2)$: $(1,2,2,3,3\dots,k-3,k-3)$ and $(1,2,2,3,3,\dots,k-4,k-3,k-3,k-2)$.
  \end{itemize}
		Thus the two vertices with degree one on layer $(s+2)$ are different from the others, which implies that $\phi$ maps the $V(G^*_{k-2})$ to itself. By Observation \ref{obs1} and  by Lemma \ref{lem6}, we obtain that $\phi$ restricted to layer $(s+2)$ is identity. 
	\end{proof}
	
	\begin{claim}
		\label{clm2}
		If the only automorphism $\phi$ of $G_{k,s}$ restricted to layer $(l+1)$, $l\ge 1$, is identity, then $\phi$ restricted to layer $l$ is also identity.
	\end{claim}
	
	\begin{proof}
		Since the automorphism $\phi$ of $G_{k,s}$ restricted to layer $(l+1)$ is identity, the corresponding edges $M$ of each edge $M_{i_l}$ in $G(i_{l+1},i_{l+2},\dots,i_{s+1})$ on layer $(l+1)$ are pairwise different. It implies that the copies $G(i_l,i_{l+1},\dots,i_{s+1})$ on layer $l$ are pairwise different and $\phi$ maps each $\widetilde{G}(i_l,i_{l+1},\dots,i_{s+1})$ to itself (if $l=1$, $\phi$ maps each $G(i_1,i_2,\dots,i_{s+1})$ to itself) and leaves the head vertex and the tail two vertices of $G(i_l,i_{l+1},\dots,i_{s+1})$ invariant. 
		By Observation \ref{obs1},and by Lemma \ref{lem5}, the automorphism $\phi$ of $G$ restricted to each vertex subset $V(G(i_l,i_{l+1},\dots,i_{s+1}))$ 
		on layer $l$ is identity. 
	\end{proof}
	
	Claim \ref{clm1} states that the automorphism $\phi$ of $G_{k,s}$ induced on layer $(s+2)$ is identity. Then by Claim \ref{clm2}, $\phi$ of $G_{k,s}$ restricted to layer $(s+1)$ is identity. Continuing this way, we obtain that $\phi$ restricted to layer $i$ is identity, $i\in\{1,2,\dots,s+2\}$. Thus $G_{k,s}$ is asymmetric.
	
	\bigskip
	
	The involution property of Theorem \ref{main4}, follows from the following claim.
	
	\begin{claim}
		\label{clm3}
		For every $k\ge 6$ and $s\ge 1$, any proper sub-$k$-graph of $G_{k,s}$ with at least $2$ vertices has an involution.
	\end{claim}
	
	\begin{proof}
		For contradiction, assume that $G_{k,s}$ contains a non-trivial sub-$k$-graph $H$ such that $H$ has no involution. Without loss of generality, let us assume that $H$ is connected.
		
		Let $l$ be the minimal layer such that there exists a copy $G=G(i_l,i_{l+1},\dots,i_{s+1})$ of $G_{k-2}$ ($G$ is a copy of $G_k$ if $l=1$ and $G=G^*_{k-2}$ if $l=s+2$) with $1<|V(H)\cap V(G)|<|V(G)|$. Let $G'$ be the sub-$(k-2)$-graph (or sub-$k$-graph) of $G$ induced by $V(G')=V(H)\cap V(G)$, and let $\widetilde{G'}$ be the corresponding sub-$k$-graph of $G'$ in $H$.
       We distinguish two cases.
		
		\noindent
		{\bf Case 1.} Such an $l$ exists.
		
		\noindent
		Let $x$, $y$ be the tail two vertices of $G$.
		If $G'$ is an empty graph, then $V(G')=\{x,y\}$, hence $H$ has an involution interchanging $x$ and $y$. Assume that $G'$ is a non-trivial sub-$(k-2)$-graph (or sub-$k$-graph) of $G$. If $G\ne G^*_{k-2}$, by Lemma \ref{lem3}, $G'$ has an involution which leaves $x$, $y$ invariant if $x$ or $y$ belongs to $V(G')$.  If $G=G^*_{k-2}$, by Lemma \ref{lem4}, $G'$ has an involution. Then by Observation \ref{obs1}, $\widetilde{G'}$ has an involution $\phi$ that maps $V(G')$ to $V(G')$, which can be easily extended to $H$.

		\noindent
		{\bf Case 2.} Such an $l$ does not exist. This means for each copy $G=G(i_m,i_{m+1},\dots,i_{s+1})$, $i\in \{1,2,\dots,s+1\}$, or $G=G^*_{k-2}$, either $V(G)\cap V(H)=V(G)$ and $\widetilde{G}$ is a sub-$(k-2)$-graph (or sub-$k$-graph) of $H$ or $V(G)\cap V(H)=\emptyset$. 
		
		\noindent
		Assume $p$ is the maximal layer such that there is a copy $G'=G(i_p,i_{p+1},\dots,i_{s+1})$ of $G_{k-2}$ ($G'$ is a copy of $G_k$ if $p=1$ and $G'=G^*_{k-2}$ if $p=s+2$) with $V(G')\cap V(H)=V(G')$. It is easy to check that the vertices of every $G(i_q, i_{q+1},\dots,i_{p-1},i_p,\dots,i_{s+1})$ ($i_q\in\{1,2,\dots,k-2\}$, $2\le q\le p-1$ and $i_q\in\{1,2,\dots, k-1\}$ if $q=1$) is contained in $H$, otherwise there exists a copy $G(i_t, i_{t+1},\dots,i_{p-1},i_p,\dots,i_{s+1})$ for some $2\le t\le p-1$, the vertices of which contained in $H$ are the tail two vertices, a contradiction. 
		Since $H$ is a non-trivial sub-$k$-graph of $G_{k,s}$, $G'\ne G^*_{k-2}$. 
		By Lemma \ref{lem3}, $G'$ has an involution. Then by Observation \ref{obs1}, $\widetilde{G'}$ has an involution $\phi$ that maps $V(G')$ to $V(G')$, which can be extended to $H$. 
	\end{proof}	
	
	\medskip
	This concludes the proof of Theorem \ref{main4}.


	\section{Concluding remarks}
	\label{sec:open}
	
	\noindent
	{\bfseries{1.}}
	Of course one can define the notion of asymmetric graph also for directed graphs.
	
	One has then the following analogy of Theorem 1:
	there are exactly $19$ minimal asymmetric binary relations.
	(These are symmetric orientations of $18$ minimal asymmetric (undirected) graphs and the single arc graph ($\{0,1\}$, $\{(0,1)\}$).)
	
	Here is a companion problem about extremal asymmetric oriented graphs. This is one of the original motivation, see e.g. \cite{BOOK}.
	
	Let $G=(V,E)$ be an asymmetric graph with at least two vertices. We say that $G$ is \emph{critical asymmetric} if for every $x\in V$ the graph $G-x=(V\setminus\{x\},\{e\in E;x\notin e\})$ fails to be asymmetric or it is exactly a single vertex.
	Recall that an oriented graph is a relation not containing two opposite arcs.
	
	\smallskip
	\noindent
	{\bfseries{Conjecture 1}}
	There is no critical oriented asymmetric graph.
 
 Explicitly: For every oriented asymmetric graph $G$ with at least two vertices, there exists $x\in V(G)$ such that $G-x$ is asymmetric. 
	
	\smallskip
	W\'{o}jcik \cite{Woj96} proved that a critical oriented asymmetric graph has to contain a directed cycle. In general, Conjecture 1 is open.
	
	\bigskip
	
	\noindent
	{\bfseries{2.}}
	More generally, we could consider $k$-ary relational structures $(X,R)$. We say the multiplicity $m(R)$ of a relation $R$ is at most $s$ if on every $k$-set there are at most $s$ tuples, $1\le s\le k!$. Thus oriented graphs are binary relations with multiplicity $1$.
    It is natural to ask for which multiplicities there are 	
 finitely many minimal asymmetric $k$-ary relational structures $(X,R)$.

	There is exactly one minimal asymmetric $k$-ary relational structure with multiplicity $1$, which is a single $k$-set. And by Theorem \ref{main3} we know there are infinitely many minimal asymmetric $k$-ary relational structures $(X,R)$ with $m(R)=k!$.
	
	Note that if a $k$-ary relation $R$ has multiplicity $m(R)=2$, then on every $k$-set the $2$ tuples should be different exactly at two places (if on every $k$-set the $2$ tuples are not in this form, then $R$ restricted to these $k$-tuples are asymmetric, which means that the minimal asymmetric $k$-ary relational structure with such $R$ is a single $k$-set). For example, for every $k$-set $\{x_1,x_2,x_3,\dots,x_k\}$, $(x_1,x_2,x_3,\dots,x_k)\in R$ implies $(x_2,x_1,x_3,\dots,x_k)\in R$. 
	
 We use the construction $G^{\circ}_{3,t}$ in Section 3.1 to prove that there are infinitely many minimal asymmetric $k$-ary relational structures ($k\ge 3$) with multiplicity $2$. Since in the proof it makes no diference if the two tuples of $R$ differ at different places, we assume that in $R$ for every $k$-set the two tuples differ at the first two places.
	
	We first construct infinitely many minimal asymmetric ternary relational structures $(X_{3,t}, R'_{3,t})$ such that for every $3$-set $\{x_1,x_2,x_3\}$, $(x_1,x_2,x_3)\in R'_{3,t}$ implies $(x_2,x_1,x_3)\in R'_{3,t}$. We use the construction of $G^{\circ}_{3,t}$ as before in Section \ref{sec3}.
	
	$G_{3,t}=(X_{3,t},\mathscr{E}_{3,t})$, 
	
	
	
	$G^{\circ}_{3,t}=((X'_{3,t},\mathscr{E}'_{3,t}))=(X_{3,t}\cup\{x\}, \mathscr{E}_{3,t}\cup\{E_{3t}\})$, where $E_{3t}=\{v_0,u_0,x\}$. 

 For every set $\{u,v,w\}\in\mathscr{E}'_{3,t}$, we have $(u,v,w)\in  R'_{3,t}$ and $(v,u,w)\in  R'_{3,t}$.
	
	The proof that $G^{\circ}_{3,t}$ is minimal asymmetric in Section 3 also works here. And then we obtain infinitely many minimal asymmetric $k$-ary relational structures $(X, R)$ with $m(R)=2$ by adding the extra $k-3$ (if $k>3$) vertices separately to each corresponding hyperedge as follows.
	
	$H^{\circ}_{k,t}=(X'_{k,t}, \mathscr{M}_{k,t})$, where $X'_{k,t}=X_{3,t}\cup\{x\}\cup\bigcup\limits_{i=0}^{3t}\{w^1_i, w^2_i,\dots,w^{k-3}_i\}$ and $ \mathscr{M}_{k,t}=\{E'_i=E_i\cup\{w^1_i, w^2_i,\dots,w^{k-3}_i\}; i\in[3t+1]\}$.
	
	Every vertex $w^1_i, w^2_i,\dots,w^{k-3}_i$ in $E_i$ for every $i\in[3t+1]$ maps to itself in any automorphism of $H^{\circ}_{k,t}$ according to the multiplicity of $R$, which complete the proof.
	
This also implies that for $k$-ary relation $R$ with $m(R)=k!-1$ there is only one minimal asymmetric relation while for $m(R)=k!-2$ we have infinitely many of them. Perhaps for every $m(R)$, $2\le m(R)\le K!-2$ there are infinitely many minimal asymmetric relations. 

Of interest are special cases such as cyclic relations.
 We call a relation $R$ {\em cyclic} if it has multiplicity $k$ and on every $k$-set $\{x_1,x_2,\dots, x_k\}$ it contains all the following tuples $(x_1,x_2,x_3,\dots,x_k)$, $(x_2,x_3,\dots,x_k,x_1)$, $\cdots$, $(x_k,x_1,x_2,\dots,x_{k-1})$. 
	
	\smallskip
	\noindent
	{\bfseries{Problem 2}}
	Are there finitely many minimal asymmetric $k$-ary cyclic relational structures $(X,R)$?
	
	It is not clear even for $k=3$.
	
\bigskip
\noindent
{\bfseries{Acknowledgement.}} The authors thank Dominik Bohnert and Christian Winter for finding a mistake in the original statement of Lemma 9.

\end{document}